\documentstyle[11pt,a4wide]{article}
\font \tenmsb=msbm10 scaled \magstep 1
\font \sevenmsb=msbm7 scaled \magstep 1
\font \fivemsb=msbm5 scaled \magstep 1
\newfam \msbfam
\textfont \msbfam=\tenmsb
\scriptfont \msbfam=\sevenmsb
\scriptscriptfont \msbfam=\fivemsb
\def \Bbb#1{\fam \msbfam \relax#1}

\font\teneufm=eufm10 scaled \magstep 1
\font\seveneufm=eufm7 scaled \magstep 1
\font\fiveeufm=eufm5 scaled \magstep 1
\newfam\eufmfam
\textfont\eufmfam=\teneufm
\scriptfont\eufmfam=\seveneufm
\scriptscriptfont\eufmfam=\fiveeufm
\def\frak#1{{\fam\eufmfam\relax#1}}

\title{\bf ON FUNCTION THEORY IN QUANTUM DISC: COVARIANCE}

\author{\sl D. Shklyarov \and \sl S. Sinel'shchikov
\and \sl L. Vaksman \thanks{Partially supported by the grant
INTAS-94-4720}}

\date{\tt Institute for Low Temperature Physics \& Engineering\\
National Academy of Sciences of Ukraine}

\newtheorem{theorem}{Theorem}[section]
\newtheorem{lemma}[theorem]{Lemma}
\newtheorem{proposition}[theorem]{Proposition}
\newtheorem{corollary}[theorem]{Corollary}

\begin{document}

\maketitle

\bigskip

\section{Modules over a Hopf algebra}

  We follow the approach of V. G. Drinfeld and M. Jimbo to constructing the
quantum groups theory \cite{J}. Everywhere in the sequel the deformation
parameter $q$ will be assumed to be a number from the interval $(0,1)$.

 The quantum universal enveloping algebra $U_q \frak{sl}_2$ is a Hopf
algebra over ${\Bbb C}$ determined by the generators $K,K^{-1}, E, F$ and
the relations
$$KK^{-1}=K^{-1}K=1,\quad K^{\pm 1}E=q^{\pm 2}EK^{\pm 1},\quad K^{\pm
1}F=q^{\mp 2}FK^{\pm 1},$$
$$EF-FE=(K-K^{-1})/(q-q^{-1}),$$
$$\Delta(K^{\pm 1})=K^{\pm 1}\otimes K^{\pm 1},\quad \Delta(E)=E \otimes 1+K
\otimes E,\quad \Delta(F)=F \otimes K^{-1}+1 \otimes F.$$

 Note that
$$\varepsilon(E)=\varepsilon(F)=\varepsilon(K^{\pm 1}-1)=0,$$
$$S(K^{\pm 1})=K^{\mp 1},\quad S(E)=-K^{-1}E,\quad S(F)=-FK,$$
with $\varepsilon:U_q \frak{sl}_2 \to{\Bbb C}$ and $S:U_q \frak{sl}_2 \to
U_q \frak{sl}_2$ being respectively the counit and the antipode of
$U_q \frak{sl}_2$.

 It was shown in \cite{CP} that $U_q \frak{sl}_2$ can be derived from the
(topological) Hopf algebra $U_h \frak{sl}_2$ over the ring of formal series
${\Bbb C}[[h]]$. The latter Hopf algebra is determined by its generators
$X^+$, $X^-$, $H$ and the relations
$$[H,X^\pm]=\pm 2X^\pm,\quad[X^+, X^-]={\rm sh}(hH/2)/{\rm sh}(h/2),$$
$$\Delta(H)=H \otimes 1+1 \otimes H,\quad \Delta(X^\pm)=X^\pm \otimes
e^{hH/4}+e^{-hH/4}\otimes X^\pm,$$
$$\varepsilon(H)=\varepsilon(X^\pm)=0,\quad S(H)=-H,\quad S(X^\pm)=-e^{\pm
h/2}X^\pm.$$
 $U_h \frak{sl}_2$ is a deformation of the ordinary universal enveloping
algebra $U \frak{sl}_2$, and the formal passage to a limit as $h \to 0$ in
the determining relations of $U_h \frak{sl}_2$ yields $U \frak{sl}_2 \simeq
U_h \frak{sl}_2/h \cdot U_h \frak{sl}_2$.

 All the $U_q \frak{sl}_2$-modules in the sequel will be assumed to be
${\Bbb R}$-graded, with $K^{-1}v=\exp(\deg(v)h/2)\cdot v$ for any
homogeneous element $v$. This restriction allows one to replace $U_q
\frak{sl}_2$ by $U_h \frak{sl}_2$ while passing to the limit as $q \to 1$ as
well as in tedious calculations:
$$q=e^{-h/2},\quad K=e^{-hH/2},\quad E=X^+e^{-hH/4},\quad
F=e^{hH/4}X^-.\eqno(1.1)$$

 Everywhere below $A$ will stand either for the Hopf algebra $U_q
\frak{sl}_2$ or the subalgebra $U_q{\frak b}_+\subset U_q \frak{sl}_2$,
generated by $K$, $K^{-1}$, $E$.

 Tensor product of $A$-modules $V_1$, $V_2$ is defined as follows:
$$A \stackrel{\Delta}{\to}A \otimes A \to {\rm End}(V_1)\otimes {\rm
End}(V_2)\simeq{\rm End}(V_1 \otimes V_2),$$ and the trivial $A$-module
${\Bbb C}$ as $$A \stackrel{\varepsilon}{\to}{\rm End}({\Bbb C})\simeq{\Bbb
C}.$$

 The morphisms of $A$-modules $\eta:V \to{\Bbb C}$ are also called invariant
integrals.

 The dual $A$-module $V^*$ is defined via the antipode $S$ as follows:
$$(al)(v)=l(S(a)v),\quad \forall a \in A,\;v \in V,\;l \in V^*.$$
It follows from the definition of the antipode $S$ that the natural pairing
$V^*\otimes V \to{\Bbb C}$ is an invariant integral.

\medskip

\begin{proposition} Linear functional $\eta:V_1 \otimes V_2 \to{\Bbb C}$ is
an invariant integral iff
$$\eta(av_1 \otimes v_2)=\eta(v_1 \otimes S(a)v_2)\quad{\rm for \;all}\quad
a \in A,\;v_1 \in V_1,\;v_2 \in V_2.\eqno(1.2)$$ \end{proposition}

\smallskip

 {\bf Proof.} (1.2) means that the map $V_1 \to V_2^*$ given by the linear
functional $\eta$ is a morphism of $A$-modules. Hence it suffices to verify
(1.2) for the generators $K^{\pm 1}$, $E$, $F$, which is straightforward. On
the other hand, in the case (1.2) is satisfied, $\eta$ becomes an invariant
integral due to the definition of the antipode $S$. \hfill $\Box$

\medskip

 Given $A$-modules $V_1,V_2$, we consider the vector space of finite
dimensional linear operators $${\rm Hom}_{\Bbb C}(V_1,V_2)_f=\{L \in{\rm
Hom}_{\Bbb C}(V_1,V_2)|\,{\rm dim}(LV_1)<\infty \}.$$

\medskip

\begin{proposition} There exists a unique $A$-module structure on ${\rm
Hom}_{\Bbb C}(V_1,V_2)_f$ such that the linear map ${\rm Hom}_{\Bbb
C}(V_1,V_2)_f \otimes V_1 \to V_2$, $L \otimes v \mapsto Lv$, $L \in {\rm
Hom}_{\Bbb C}(V_1,V_2)_f$, $v \in V_1$, is an $A$-module morphism.
\end{proposition}

\smallskip

 {\bf Proof.} The existence is evident in view of the canonical isomorphism
of the vector spaces $i:{\rm Hom}_{\Bbb C}(V_1,V_2)_f {\displaystyle \to
\atop \displaystyle \sim} V_2 \otimes V_1^*$. It remains to demonstrate that
for any $A$-module structure on ${\rm Hom}_{\Bbb C}(V_1,V_2)_f$ with the
above properties, $i$ is a morphism of $A$-modules. This can be easily
verified for the generators $K$, $K^{-1}$, $E$, $F$. \hfill $\Box$

\medskip

\begin{corollary} The finite dimensional linear operators in $V$ form an
$A$-module ${\rm End}(V_f)\simeq V \otimes V^*$. \hfill $\Box$
\end{corollary}

\medskip

\begin{proposition} The linear functional
$${\rm End}(V_f)\to{\Bbb C},\quad L \mapsto {\rm tr}(LK^{-1})\eqno(1.3)$$
is an invariant integral.
\end{proposition}

\smallskip

 {\bf Proof.} Let $i_{st}:V \hookrightarrow V^{**}$ be the standard
embedding of vector spaces. It follows from the relation $S^2(a)=K^{-1}aK$
that the map $i:V \hookrightarrow V^{**}$, $i:v \mapsto i_{st}\cdot K^{-1}v$
is an $A$-module morphism. It induces an invariant integral
$${\rm End}(V)_f \simeq V \otimes V^*\stackrel{i \otimes{\rm
id}}{\hookrightarrow}V^{**}\otimes V^*\to{\Bbb C}.$$
Evidently, this integral coincides with that in (1.3). \hfill $\Box$

\medskip

 To conclude, note that the tensor product of $A$-module morphisms is again
an $A$-module morphism. More exact statement is that $A$-modules and their
morphisms constitute a tensor category \cite{CP}.

\bigskip

\section{Covariant algebras and bimodules}

 Consider an algebra $F$, which is also an $A$-module. $F$ is said to be an
$A$-module (covariant) algebra if the multiplication
$$m:F \otimes F \to F,\quad m:f_1 \otimes f_2 \mapsto f_1f_2;\quad f_1,f_2
\in F,$$
is an $A$-module morphism \cite{A}.

 In the case of a unital algebra $F$ the above definition also includes the
assumption
$$\forall a \in A \qquad a \cdot 1=\varepsilon(a)1.\eqno(2.1)$$

 An element $v$ of a $A$-module $V$ is called an invariant if the map ${\Bbb
C}\to V$, $z \mapsto z \cdot v$, $z \in{\Bbb C}$, is an $A$-module morphism.
In this context, (2.1) claims that the unit of $F$ is an invariant.

 One can view the algebra ${\rm End}(V)_f \simeq V \otimes V^*$ and the
tensor algebra $T(V)={\Bbb C}+V+V \otimes V+\ldots$,\ associated to the
$A$-module $V$, as the examples of covariant algebras.

 A bimodule $M$ over a covariant algebra $F$ is said to be $A$-module (or
covariant) if $M$ itself is an $A$-module, and the maps
$$F \otimes M \to M,\quad f \otimes m \mapsto fm,$$
$$M \otimes F \to M,\quad m \otimes f \mapsto mf,\quad f \in F,\;m \in M,$$
are $A$-module morphisms.

 Impose the notation $U_q \frak{su}(1,1)=(U_q \frak{sl}_2,*)$, $U_h
\frak{su}(1,1)=(U_h \frak{sl}_2,*)$ for the Hopf *-algebras whose
involutions are given by
$$h^*=h,\quad H^*=H,\quad (X^\pm)^*=-X^\mp,$$
$$E^*=-KF,\quad F^*=-EK^{-1},\quad K^*=K.$$

 A covariant algebra $F$ is called a covariant *-algebra if it is equipped
with an involution, and
$$\forall f \in F \qquad (\xi f)^*=(S(\xi))^*f^*\eqno(2.2)$$
for all elements $\xi$ of the *-Hopf algebra.

 It is very well known \cite{KL, NN} that ${\rm Pol}({\Bbb C})_q$ can be
equipped with a structure of a covariant *-algebra. This result was extended
onto the case of "prehomogeneous vector spaces of parabolic type" in
\cite{SV}, where the detailed calculations were also performed for the case
of ${\rm Pol}({\Bbb C})_q$. In particular, one has \cite{SV}:
$${K^\pm z=q^{\pm 2}z,\qquad Fz=q^{1/2},\qquad Ez=-q^{1/2}z^2,\atop K^{\pm
1}z^*=q^{\mp 2}z^*,\qquad Ez^*=q^{-3/2},\qquad
Fz^*=-q^{5/2}z^{*2};}\eqno(2.3)$$
("in a different notation" we have
$${Hz=2z,\qquad X^-z=e^{-h/4},\qquad X^+z=-e^{h/4}z^2,\atop Hz^*=-2z^*,\qquad
X^+z^*=e^{h/4},\qquad X^-z^*=-e^{-h/4}z^{*2}).}\eqno(2.4)$$

 Let $F$ be a covariant algebra, $M$ a covariant bimodule over $F$, and
$\eta:M \to{\Bbb C}$ --- an invariant integral. Proposition 1.1 implies the
following "formula of integration by parts":

\medskip

\begin{proposition} $\forall f \in F$, $\psi \in M$, $a \in A$
$$\int(af)\psi d \eta=\int f(S(a)\psi)d \eta,\qquad \int(a \psi)fd \eta=\int
\psi(S(a)f)d \eta.$$
\end{proposition}

\bigskip

\section{The covariant algebra $\bf D(U)_q$}

 We refer to \cite{SSV1} for definitions of ${\rm Pol}({\Bbb
C})_q$-bimodules $D(U)_q$, $D(U)_q^\prime$ of finite functions and
distributions respectively in the quantum disc. It follows from the
definitions that $D(U)_q^\prime$ may be identified with the space of formal
series
$$f=\sum_{jk}a_{jk}z^jz^{*k}.\eqno(3.1)$$
The subspace ${\rm Pol}({\Bbb C})_q$ is constituted by finite sums of the
form (3.1), and the structure of algebra in ${\rm Pol}({\Bbb C})_q$ and
that of  bimodule in $D(U)_q^\prime$ are given by the commutation relation
$z^*z=q^2zz^*+1-q^2$. It follows from the definition of $D(U)_q$ that $f \in
D(U)_q$ iff
$$z^{*N} \cdot f=f \cdot z^N=0 \eqno(3.2)$$
for some $N \in{\Bbb N}$. The space $D(U)_q^\prime$ is equipped with the
topology of coefficientwise convergence, so ${\rm Pol}({\Bbb C})_q$ and
$D(U)_q$ become its dense linear subspaces (see \cite{SSV1}).

 It follows from the covariance of ${\rm Pol}({\Bbb C})_q$, together with
(2.3), (3.1), (3.2), that the action of the Hopf algebra $U_q
\frak{sl}(1,1)$ can be extended by continuity onto $D(U)_q^\prime$ and hence
transferred onto $D(U)_q$. Of course, $D(U)_q$ is a covariant $*$-algebra,
and $D(U)_q^\prime$  a covariant $D(U)_q$-bimodule.

 The work \cite{SSV1} contains a construction of an isomorphism $T$ between
$D(U)_q$ and the algebra of matrices $(l_{mn})_{m,n \in{\Bbb Z}_+}$ with
finitely many nonzero matrix elements. Impose a finite function $f_0$ in
the quantum disc such that
$$l_{mn}(f_0)=\left \{\begin{array}{l|l}1 & m=n=0 \\
                                        0 & m^2+n^2 \ne 0
                      \end{array}\right..$$

 Let ${\Bbb C}[z]_q$ and ${\Bbb C}[z^*]_q$ stand for covariant subalgebras
of ${\rm Pol}({\Bbb C})_q$, generated respectively by $z$ and $z^*$. It
follows from the definitions that $D(U)_q$ is a covariant ${\rm Pol}({\Bbb
C})_q$-bimodule. As for $f_0$, one evidently has $f_0 \cdot f_0=f_0$, and

\medskip

\begin{proposition}\hfill
\begin{enumerate}

\item $\{f \in D(U)_q^\prime|\;z^*f=fz=0\}={\Bbb C}f_0$

\item $D(U)_q={\Bbb C}[z]_q \cdot f_0 \cdot {\Bbb C}[z^*]_q.$

\end{enumerate}
\end{proposition}

\medskip

 Let $\widetilde{H}=\{f \in D(U)_q|\;fz=0 \}$. The relation (2.3) and the
covariance of the ${\rm Pol}({\Bbb C})_q$-bimodule $D(U)_q$ imply

\medskip

\begin{proposition} $U_q {\frak b}_+\cdot \widetilde{H} \subset
\widetilde{H}$.
\end{proposition}

\medskip

 It is also easy to deduce

\medskip

\begin{proposition} $\widetilde{H}={\Bbb C}[z]_qf_0={\rm Pol}({\Bbb
C})_qf_0=D(U)_qf_0$.
\end{proposition}

\medskip

 According to propositions 3.2 and 3.3, the vector space $\widetilde{H}$ is
a $U_q{\frak b}_+$-module and a $D(U)_q$-module.  We denote the
corresponding representations of $U_q{\frak b}_+$ and $D(U)_q$ in
$\widetilde{H}$ respectively by $\widetilde{\Gamma}$ and $\widetilde{T}$.

 We also use the notation $y=1-zz^*$. The following is a straightforward
consequence of the definitions.

\medskip

\begin{proposition} $K^{\pm 1}\psi(y)=\psi(y)$ for any $\psi(y)\in
D(U)_q^\prime$.
\end{proposition}

\medskip

\begin{theorem} The linear functional
$$\eta:D(U)_q \to{\Bbb C},\quad \eta:f \mapsto{\rm
tr}(\widetilde{T}(f)\widetilde{\Gamma}(K^{-1}))\eqno(3.3)$$
is an invariant integral.
\end{theorem}

\smallskip

 {\bf Proof.} It is easy to equip $\widetilde{H}$ with a structure of
pre-Hilbert space in such a way that $\widetilde{T}(f^*)=\widetilde{T}(f)^*$
for all $f \in D(U)_q$. By a virtue of proposition 3.4 one has
$K^{-1}f_0=f_0$. This implies $K^{-1}(z^nf_0)=q^{-2n}z^nf_0$, and hence the
linear functional $\eta$ is real: $\eta(f^*)=\overline{\eta(f)}$ for all $f
\in D(U)_q$. In view of (2.2), it suffices to show the $U_q{\frak
b}_+$-invariance of this linear functional. For that, by proposition 1.4, it
suffices to demonstrate that the morphism of algebras $\widetilde{T}:D(U)_q
\to {\rm End}(\widetilde{H})_f$ is also a morphism of $U_q{\frak
b}_+$-modules. To see that, one has to apply the covariance of $D(U)_q$ and
to use the same argument as in the proof of proposition 1.2. \hfill $\Box$

\medskip

 {\bf Remark 3.6.} The non-negativity of the invariant integral $\eta$
follows from the existence of a positive scalar product under which
$\widetilde{T}(z^*)=\widetilde{T}(z)^*$ (see section 1 of \cite{SSV1}).

\medskip \stepcounter{theorem}

 Recall \cite{SSV1} that each element $f \in D(U)_q$ admits a unique
decomposition
$$f=\sum_{j>0}z^j
\psi_j(y)+\psi_0(y)+\sum_{j>0}\psi_{-j}(y)z^{*j}.\eqno(3.4)$$
Here $\psi_j$ are finite functions defined on $q^{2{\Bbb Z}_+}$ such that
${\tt card}\{(i,j)|\;\psi_j(q^{2i})\ne 0 \}<\infty$.

\medskip

 {\bf Remark 3.7.} The linear functional
$$\int \limits_{U_q}fd \nu=(1-q^2)\sum_{j=0}^\infty \psi_0(q^{2j})q^{-2j}$$
used in \cite{SSV1} differs only by a scalar multiple from the linear
functional $\eta$, and hence is an invariant integral.

\medskip \stepcounter{theorem}

\begin{proposition}
$$X^+f_0=-\frac{e^{-h/4}}{1-e^{-h}}zf_0,\qquad
X^-f_0=-\frac{e^{-3h/4}}{1-e^{-h}}f_0z^*.\eqno(3.5)$$
\end{proposition}

\smallskip

 {\bf Proof.} In virtue of (2.2), it suffices to verify the first relation.
It was shown above that $X^+f_0 \subset \widetilde{H}={\Bbb
C}[z]f_0$. Besides, $H(X^+f_0)=2X^+f_0$. This implies $X^+f_0=c_+zf_0$ for
some constant $c_+$. The value of $c_+$ can be found via an application of
$X^+$ to the relation $z^*\cdot f_0=0$:
$$0=X^+(z^*f_0)=e^{h/4}f_0+e^{h/2}z^*c_+zf_0,$$
$$0=e^{h/4}+e^{h/2}(1-e^{-h})c_+.\eqno \Box$$

\medskip

 The proofs of the basic results announced in \cite{SSV1} require the
following statement.

\medskip

\begin{theorem} $f_0$ generates the $U_q \frak{sl}_2$-module $D(U)_q$.
\end{theorem}

\smallskip

 {\bf Proof.} Use the covariance of $D(U)_q$, proposition 3.8 and (2.4) to
get
$$X^+(z^jf_0)=a_jz^{j+1}f_0,\eqno(3.6)$$
$$\left \{{X^-(f_0z^{*k})=b_kf_0z^{*(k+1)}\atop
X^-(z^{j+1}f_0z^{*k})=c_{jk}z^jf_0z^{*k}+d_{jk}z^{j+1}f_0z^{*(k+1)},}\right.
\eqno(3.7)$$
and the inequalities $a_j<0$, $b_k<0$, $c_{jk}>0$, $d_{jk}<0$ for all $j,k
\in{\Bbb Z}_+$.

 Consider the subspaces $L_m$, $m \in{\Bbb Z}_+$:
$$L_m=\{f \in D(U)_q|\;f=\sum_{j=0}^m \sum_{i=0}^\infty
a_{ij}z^if_0z^{*j},\;a_{ij}\in{\Bbb C}\}.$$
It suffices to prove that for all $m$, $U_q \frak{sl}_2 \cdot f_0 \supset
L_m$. We proceed by induction. Firstly, $L_0=U_q \frak{b}_+\cdot f_0$
because of (3.6) and $a_j \ne 0$ for all $j \in{\Bbb Z}_+$. Secondly,
$L_{m+1}\subset L_m+X^-L_m$ due to (3.7) and $b_k \ne 0$, $c_{jk}\ne 0$,
$d_{jk}\ne 0$ for all $j,k \in{\Bbb Z}_+$. \hfill $\Box$

\medskip

\begin{corollary} (Uniqueness of the invariant integral)

${\rm dim \;Hom}_{U_q \frak{sl}_2}(D(U)_q,{\Bbb C}) \le 1$.
\end{corollary}

\smallskip

 {\bf Proof.} If $V$ is a free $U_q \frak{sl}_2$-module, ${\rm
dim \;Hom}_{U_q \frak{sl}_2}(V,{\Bbb C})=1$. On the other hand, theorem 3.9
assures the embedding of the vector spaces ${\rm Hom}_{U_q
\frak{sl}_2}(D(U)_q,{\Bbb C})\hookrightarrow {\rm Hom}_{U_q
\frak{sl}_2}(V,{\Bbb C})$. \hfill $\Box$

\medskip

 Finally, note that one can prove the existence and the uniqueness of an
invariant integral (${\rm dim \;Hom}_{U_q \frak{sl}_2}(D(U)_q,{\Bbb C})=1$)
in many different ways. These are not those facts themselves that make an
interest, but the explicit formula (3.3) for an invariant integral. A
similar formula was published in \cite{KV}.

\bigskip

\section{Invariance of the operator $\Box$}

 Consider a module $V$ over a Hopf *-algebra $A$. By its definition, the
antimodule $\overline{V}$ coincides with $V$ as an Abelian group, but the
multiplication by scalars and the $A$-action in $\overline{V}$ are given by
$$(\lambda,\overline{v}) \mapsto \overline{\lambda}\cdot
\overline{v};\quad(a,\overline{v})\mapsto(S(a))^*\overline{v}.$$

 A sesquilinear form $V_1 \times V_2 \to{\Bbb C}$ is said to be invariant if
the associated linear functional $\eta:\overline{V}_2 \otimes V_1 \to{\Bbb
C}$ is an invariant integral.

 It follows from the invariance of the scalar product in $V$ that
$$\forall a \in A,\;v^\prime,v^{\prime \prime}\in V \qquad
(av^\prime,v^{\prime \prime})=(v^\prime,a^*v^{\prime \prime}).$$
In fact, by a virtue of proposition 1.1,
$$(av^\prime,v^{\prime \prime})=\eta(\overline{v^{\prime \prime}}\otimes
av^\prime)=\eta(S^{-1}(a)\overline{v^{\prime \prime}}\otimes
v^\prime)=\eta(\overline{a^*v^{\prime \prime}}\otimes
v^\prime)=(v^\prime,a^*v^{\prime \prime}).$$

\medskip

\begin{proposition} The scalar product $(f_1,f_2)=\displaystyle \int
\limits_{U_q}f_2^*f_1d \nu$ in $D(U)_q$ is invariant.
\end{proposition}

\smallskip

 {\bf Proof.} This is a straightforward consequence of invariance of the
integral $\displaystyle \int \limits_{U_q}fd \nu$ proved in the previous
section. \hfill $\Box$

\medskip

 It is easy to show that the structure of a covariant algebra in
$\Omega({\Bbb C})_q$ introduced in \cite{SV} can be transferred by a
continuity onto $\Omega(U)_q$. The differential $d:\Omega(U)_q \to
\Omega(U)_q$ is certainly a morphism of $U_q \frak{sl}_2$-modules.

 Consider the integral $\Omega(U)_q^{(1,1)}\to{\Bbb C}$, $\displaystyle \int
\limits_{U_q}fdz^*dz \stackrel{def}{=}q^2 \cdot2i\pi \displaystyle \int
\limits_{U_q}f(1-zz^*)^2d \nu$ introduced in \cite[section 4]{SSV1}. This
linear functional is an invariant integral since the linear operator
$$\Omega(U)_q^{(1,1)}\to D(U)_q;\quad fdz^*dz \mapsto f \cdot(1-zz^*)^2$$
is a morphism of $U_q \frak{sl}_2$-modules.

 The work \cite{SSV1} implements the bimodules $\Omega({\Bbb
C})_{\lambda,q}^{(0,*)}$, $\lambda \in{\Bbb R}$, over the algebra
$\Omega({\Bbb C})_q^{(0,*)}$. It was shown in \cite{SV} that these
bimodules are covariant, and the operator $\overline{\partial}:\Omega({\Bbb
C})_{\lambda,q}^{(0,*)}\to \Omega({\Bbb C})_{\lambda,q}^{(0,*)}$ is a
morphism of $U_q \frak{sl}_2$-modules. Furthermore, for the generator $v_
\lambda$ involved into the definition of $\Omega({\Bbb
C})_{\lambda,q}^{(0,*)}$, one has (\cite{SV}): $\overline{\partial}v_
\lambda=0$,
$$ K^{\pm 1}v_ \lambda=q^{\pm \lambda}\cdot v_ \lambda,\quad Fv_
\lambda=0.$$
(One can deduce also $Ev_ \lambda=-q^{1/2}\frac{\textstyle 1-q^{2
\lambda}}{\textstyle 1-q^2}zv_ \lambda$ via observing that $Ev_ \lambda \ne
0$, $Ev_ \lambda={\rm const}\cdot zv_ \lambda$, $F(Ev_ \lambda)=-(EF-FE)v_
\lambda=-\frac{\textstyle q^\lambda-q^{-\lambda}}{\textstyle q-q^{-1}}v_
\lambda$, \ \ $F(zv_ \lambda)=F(z)(K^{-1}v_ \lambda)=q^{1/2}q^{-\lambda}v_
\lambda$.)

 One can find in \cite{SSV1} the construction of the covariant algebra
$\Omega(U)_q^{(0,*)}$ and the bimodules $\Omega(U)_{\lambda,q}^{(0,*)}$
over this algebra as completions in some special topology. It is easy to
show that the actions of the Hopf algebra $U_q \frak{sl}_2$ and the operator
$\overline{\partial}$ can be transferred by a continuity from $\Omega({\Bbb
C})_{\lambda,q}^{(0,*)}$ onto $\Omega(U)_{\lambda,q}^{(0,*)}$. Thus we
obtain the covariant algebra, covariant bimodules and morphisms
$\overline{\partial}:\Omega(U)_{\lambda,q}^{(0,*)}\to
\Omega(U)_{\lambda,q}^{(0,*)}$ of $U_q \frak{sl}_2$-modules.

 Concerning the following proposition, we refer the reader to the relations
(2.4), (2.5) of \cite{SSV1} which determine the scalar products in
$\Omega(U)_{\lambda,q}^{(0,0)}$ and $\Omega(U)_{\lambda,q}^{(0,1)}$.

\medskip

\begin{proposition} The scalar products in $\Omega(U)_{\lambda,q}^{(0,0)}$
and $\Omega(U)_{\lambda,q}^{(0,1)}$ are invariant.
\end{proposition}

\smallskip

 {\bf Proof.} One can show that the linear operators
$$j_0:\overline{\Omega(U)_{\lambda,q}^{(0,0)}}\otimes
\Omega(U)_{\lambda,q}^{(0,0)}\to \Omega(U)_q^{(0,0)},$$
$$j_0:\overline{(f_2 \cdot v_ \lambda)}\otimes(f_1 \cdot v_ \lambda)\mapsto
f_2^*\cdot f_1(1-zz^*)^{-\lambda},$$
$$j_1:\overline{\Omega(U)_{\lambda,q}^{(0,1)}}\otimes
\Omega(U)_{\lambda,q}^{(0,1)}\to \Omega(U)_q^{(1,1)},$$
$$j_1:\overline{(f_2v_ \lambda dz^*)}\otimes(f_1v_ \lambda dz^*)\mapsto
f_2^*\cdot f_1(1-zz^*)^{-\lambda}dzdz^*$$
are the morphisms of $U_q \frak{sl}_2$-modules. So it remains to apply the
invariance of the integrals
$$\Omega(U)_q^{(0,0)}\to{\Bbb C},\quad \Omega(U)_q^{(1,1)}\to{\Bbb C},$$
$$f \mapsto \int \limits_{U_q}fd \nu,\quad f \cdot dz^*dz \mapsto \int
\limits_{U_q}f \cdot(1-zz^*)^2d \nu,$$
which was proved above. \hfill $\Box$

\medskip

 An immediate consequence of proposition 4.2 is

\medskip

\begin{theorem} The restrictions of linear operators
$$\overline{\partial}^*:L^2(d \mu)_q \to L^2(d \nu)_q,\quad
\Box=-\overline{\partial}^*\overline{\partial}:L^2(d \nu)_q \to L^2(d
\nu)_q$$
onto the dense linear subspaces $\Omega(U)_q^{(0,1)}$ and $D(U)_q$ are
morphisms of $U_q \frak{sl}_2$-modules.
\end{theorem}

\medskip

 If $V$ is a module over the Hopf algebra $U_q \frak{sl}_2$, then the
invariance of $v \in V$ means that $Ev=Fv=(K^{\pm 1}-1)v=0$.

\medskip

\begin{proposition} The map which takes a distribution $k \in D(U)_q^\prime$
to the linear functional $D(U)_q \to{\Bbb C}$, $f \mapsto \displaystyle \int
\limits_{U_q}k \cdot fd \nu$, realizes an isomorphism between the $U_q
\frak{sl}_2$-module $D(U)_q^\prime$ and the $U_q \frak{sl}_2$-module dual to
$D(U)_q$.
\end{proposition}

\smallskip

 {\bf Proof.} It suffices to use proposition 2.1, the invariance of
integral, and the definition of a dual $U_q \frak{sl}_2$-module. \hfill
$\Box$

\medskip

 Transfer by a continuity the $U_q \frak{sl}_2$-action from the algebra
${\rm Pol}({\Bbb C})_q^{op}\otimes{\rm Pol}({\Bbb C})_q$ onto its completion
$D(U \times U)_q^\prime$.

\medskip

\begin{proposition}\hfill
\begin{enumerate}

\item The map that takes a "kernel" $K \in D(U \times U)_q^\prime$ to the
linear operator $D(U)_q \to D(U)_q^\prime$, $f \mapsto \displaystyle \int
\limits_{U_q}K(z,\zeta)f(\zeta)d \nu$ is one-to-one.

\item The integral operator $D(U)_q \to D(U)_q^\prime$ is a morphism of $U_q
\frak{sl}_2$-modules iff its kernel is invariant.
\end{enumerate}
\end{proposition}

\smallskip

 {\bf Proof.} The first statement is a direct consequence of the
definitions, and the second one follows from proposition 2.1. In fact, the
invariance of a kernel $K$ is equivalent to being a solution of the "partial
differential equation":
$$(\xi \otimes 1)K=(1 \otimes S^{-1}(\xi))K, \quad \xi \in U_q
\frak{sl}_2.\eqno \Box$$

\medskip

 {\bf Remark 4.6.} The invariance of $\Box$ and proposition 4.5 hint that
$\Box^{-1}$ is an integral operator with invariant kernel.

\bigskip

\section{The operator $\Box$ and the Casimir element $\Omega$}

 Consider the element
$$\Omega=FE+{1 \over(q^{-1}-q)^2}(q^{-1}K^{-1}+qK-(q^{-1}+q))\eqno(5.1)$$
of $U_q \frak{su}(1,1)$. In view of (1.1) one also has
$$\Omega=X^-X^++{1 \over {\rm sh}^2(h/2)}{\rm sh}(Hh/4){\rm
sh}((H+2)h/4).\eqno(5.2)$$
(At the limit $h \to 0$ we get $\Omega=X^-X^++{\textstyle 1 \over
\textstyle 4}H(H+2)$.) It is well known (see \cite{CP}) and is easy to
verify that $\Omega$ is in the center of $U_q \frak{su}(1,1)$.

 It is worthwhile to note that
$$\Omega^*=\Omega,\qquad S(\Omega)=\Omega.\eqno(5.3)$$
In fact,
$$(FE)^*=FE,\qquad(K^{\pm 1})^*=K^{\pm 1},\qquad S(FE)=EF,$$
$$S(\Omega)-\Omega={K^{-1}-K \over
q^{-1}-q}-{q^{-1}K^{-1}+qK-q^{-1}K-qK^{-1}\over(q^{-1}-q)^2}=0.$$

 Remind that the Laplace-Beltrami operator $\Box$ was defined in \cite{SSV1}
as $\Box=-\overline{\partial}^*\overline{\partial}$.

\medskip

\begin{lemma} For any function $\psi(t)$ with a finite carrier inside
$q^{-2{\Bbb Z}_+}$,
$$\Box \psi(x)=-Dx(1-q^{-1}x)D \psi(x),$$
with $x=(1-zz^*)^{-1}$, $D:f(t)\mapsto(f(q^{-1}t)-f(qt))/(q^{-1}t-qt)$.
\end{lemma}

\smallskip

 {\bf Proof.} We use here the scalar product considered in section 4 with
$\lambda=0$. To prove the lemma, it suffices to obtain the relations
$$\|\psi(x)\|^2=q^2 \int \limits_1^\infty|\psi(t)|^2d_{q^{-2}}t,\eqno(5.4)$$
$$\|\overline{\partial}\psi(x)\|^2=q^2 \int \limits_1^\infty(Dt(q^{-1}t-1)D
\psi(t))\overline{\psi(t)}d_{q^{-2}}t,\eqno(5.5)$$
with
$$\int \limits_1^\infty f(t)d_{q^{-2}}t \stackrel{\rm
def}{=}(q^{-2}-1)\sum_{m=0}^\infty f(q^{-2m})q^{-2m},$$
and $\psi$ is a function with compact support ${\tt supp}\,\psi \in q^{-2{\Bbb
Z}_+}$.

 (5.4) is obvious, and (5.5) follows from the relation
$$\overline{\partial}f(y)=-z \frac{f(y)-f(q^2y)}{y-q^2y}dz^*,\eqno(5.6)$$
with $y=1-zz^*$. (5.6) implies that
$$\overline{\partial}\psi(x)=-{z \over 1-q^2}x(\psi(x)-\psi(q^{-2}x))dz^*$$
for any function $\psi$ with a compact carrier ${\tt supp}\,\psi \in
q^{-2{\Bbb Z}_+}$. Hence
$$\|\overline{\partial}\psi(x)\|^2=(1-q^2)^{-1}\sum_{m=0}^\infty
\left|\psi(q^{-2m})-\psi(q^{-(2m+2)})\right|^2(1-q^{2m+2})q^{-2m}=$$
$$=q^{-2}(q^{-2}-1)^{-2}\int \limits_1^\infty
\left|\psi(t)-\psi(q^{-2}t)\right|^2(1-q^2t^{-1})d_{q^{-2}}t=$$
$$=-\int \limits_1^\infty
\left|\frac{\psi(t)-\psi(q^{-2}t)}{t-q^{-2}t}\right|^2 \cdot
t(1-q^{-2}t)d_{q^{-2}}t.$$

 We are to apply the "integration by parts" formula
$$\int \limits_0^\infty u_1(x)\cdot
\frac{u_2(q^{-2}x)-u_2(x)}{q^{-2}x-x}d_{q^2}x=-q^2 \int \limits_0^\infty
\frac{u_1(x)-u_1(q^2x)}{x-q^2x}\cdot u_2(x)d_{q^2}x,$$
with $\int \limits_0^\infty u(x)d_{q^2}x=(1-q^2)\sum
\limits_{m=-\infty}^\infty u(q^{2m})q^{2m}$. In this way we obtain
$$\Box \psi=-B_+t(1-q^{-2}t)B_-\psi,$$
with $B_+$, $B_-$ being the linear operators given by $B_ \pm
f(t)=\frac{\textstyle f(t)-f(q^{\pm 2}t)}{\textstyle t-q^{\pm 2}t}$.
Finally we have
$$\Box \psi=-\left.D(qt)(1-q^{-1}t)(B_-\psi)\right|_{qt}=-Dt(1-q^{-1}t)D
\psi.\eqno \Box$$

\medskip

\begin{proposition} For all $f \in D(U)_q$ one has $q \Box f=\Omega f$.
\end{proposition}

\smallskip

 {\bf Proof.} Let $f_j(x)=\left \{\begin{array}{cc}1, & x=q^{-2j} \\
                                                   0, & x \ne q^{-2j}
                                  \end{array}\right.,$\\
with $x=(1-zz^*)^{-1}$. By a virtue of theorems 3.9, 4.3, it suffices to
prove that $\Box f_0=\Omega f_0$. Now (2.4), (3.5) imply that
$$\Omega f_0=X^-X^+f_0=X^-\left(-{e^{-h/4}\over 1-e^{-h}}zf_0 \right)=$$
$$=-{e^{-h/2}\over 1-e^{-h}}f_0+{e^{-3h/2}\over(1-e^{-h})^2}zf_0z^*=-{1
\over q^{-1}-q}f_0+{q^2 \over q^{-1}-q}f_1.$$

 It remains to apply lemma 5.1:
$$\Box f_0=-Dx(1-q^{-1}x)Df_0=-{f_0 \over 1-q^2}+{q^2f_1 \over 1-q^2}.\eqno
\Box$$

\medskip

\begin{corollary} The Laplace-Beltrami operator $\Box$ is extendable by a
continuity from the linear subspace $D(U)_q \subset D(U)_q^\prime$ onto the
entire distribution space $D(U)_q^\prime$.
\end{corollary}

\medskip

\begin{corollary} For any function $\psi(t)$ on $q^{-2{\Bbb Z}_+}$ one has
$$\Omega \psi(x)=q \Box \psi(x)=qDx(q^{-1}x-1)D \psi(x),\eqno(5.7)$$
with $x=(1-zz^*)^{-1}$.
\end{corollary}

\medskip

 Consider the Hilbert space $L^2(d \nu)_q^{(0)}$ of such functions on
$q^{-2{\Bbb Z}_+}$ that $\int \limits_1^\infty
|f(x)|^2d_{q^{-2}}x<\infty$. Let $\Box^{(0)}:L^2(d \nu)_q^{(0)}\to L^2(d
\nu)_q^{(0)}$, $\Box^{(0)}:f(x)\mapsto-Dx(1-q^{-1}x)Df(x)$, be "the radial
part" of $\Box$.

\medskip

\begin{lemma} There exist constants $0<c_1 \le c_2$ such that $c_1
\le-\Box^{(0)}\le c_2$.
\end{lemma}

\smallskip

 {\bf Proof.} It follows from the definitions that the operator $\Box^{(0)}$
is selfadjoint, bounded, and non-positive: $0 \le-\Box^{(0)}\le c_2$. It
remains to show that $0$ is not in the spectrum of $\Box^{(0)}$.

 A similar result for the operator
$$q^{-1}Dx^2D:f(x)\mapsto{1 \over(q^{-1}-q)^2}\left \{
\begin{array}{cl}q^{-2}f(q^{-2}x)-(1+q^{-2})f(x)+f(q^2x), & x
                 \in q^{-2{\Bbb N}}\\
                 q^{-2}f(q^{-2})-(1+q^{-2})f(1), & x=1
\end{array}\right.$$
in $L^2(d \nu)_q^{(0)}$ is easily available via considering an
expansion in eigenfunctions. On the other hand, the continuous spectra
of the operators $\Box^{(0)}$ and $q^{-1}Dx^2D$ coincide since their
difference is a compact operator (see \cite{K}). Hence $0$ does not belong
to the continuous spectrum of $\Box^{(0)}$. It remains to prove that $0$ is
also not an eigenvalue of this operator. Otherwise, there should exist such
a function $\psi(y)\in L^2(d \nu)_q^{(0)}$ that $\overline{\partial}\psi=0$.
That is,
$$-z \frac{\psi(y)-\psi(q^2y)}{y-q^2y}dz^*=0.$$
Hence $\psi(y)$ is a constant, and $\int \limits_{U_q}|\psi|^2d \nu=\infty$.
Thus, we get a contradiction due to our assumption that $0$ is an eigenvalue
of $\Box^{(0)}$. \hfill $\Box$

\medskip

 To finish our observations, we sketch the proof of \cite[lemma 3.1]{SSV1}.
Let us prove that $c_1 \le-\Box \le c_2$.

 Consider a vector space with a basis $\{e_k \}_{k \in{\Bbb Z}}$ and a
complex number $l$. Let $V^{(l)}$ stand for a $U_q \frak{sl}_2$-module given
by
$$X^\pm e_k=\pm \frac{{\rm sh}((k \mp l)h/2)}{{\rm sh}(h/2)}e_{k \pm
1},\qquad He_k=2ke_k.$$
This parametrization of $U_q \frak{sl}_2$-modules is justified by the
possibility of realizing them in the spaces of functions with fixed
homogeneity degree $l$ on a quantum cone (see \cite{VK, SSV3}).

 It follows from the definitions that
$$\Omega v=\frac{{\rm sh}(lh/2){\rm sh}((l+1)h/2)}{{\rm sh}^2(h/2)}v,\qquad
v \in V^{(l)},\eqno(5.8)$$
and with $l \notin{\Bbb Z}$,\ $V^{(l)}$ is a simple $U_q \frak{sl}_2$-module.
We need only those $U_q \frak{sl}_2$-modules $V^{(l)}$ for which ${\rm
sh}(lh/2){\rm sh}((l+1)h/2)<0$. If $l_1-l_2 \in{\textstyle 2 \pi i \over
\textstyle h}{\Bbb Z}$ or $l_1+l_2+1 \in{\textstyle 2 \pi i \over \textstyle
h}{\Bbb Z}$, then obviously $V^{(l_1)}\approx V^{(l_2)}$.  Conversely, if
$V^{(l_1)}\approx V^{(l_2)}$, then ${\rm sh}(l_1h/2){\rm
sh}((l_1+1)h/2)={\rm sh}(l_2h/2){\rm sh}((l_2+1)h/2)$ and hence $l_1-l_2
\in{\textstyle 2 \pi i \over \textstyle h}{\Bbb Z}$ or $l_1+l_2+1
\in{\textstyle 2 \pi i \over \textstyle h}{\Bbb Z}$.

 The above observations allow one to restrict oneself to the $U_q
\frak{sl}_2$-modules $V^{(l)}$ with $l \in{\cal L}={\cal L}_1 \cup{\cal
L}_2 \cup{\cal L}_3$,
$${\cal L}_1=\{l \in{\Bbb C}|\;-{1 \over 2}<{\rm Re}\,l<0,\;{\rm
Im}\,l=0\},$$
$${\cal L}_2=\{l \in{\Bbb C}|\;{\rm Re}\,l=-{1 \over 2},\;0 \le{\rm
Im}\,l<{\pi i \over h}\},$$
$${\cal L}_3=\{l \in{\Bbb C}|\;{\rm Re}\,l>-{1 \over 2},\;{\rm Im}\,l={\pi i
\over h}\}.$$

 It is easy to prove the existence and the uniqueness of an invariant scalar
product in $V^{(l)}$ with $(e_0,e_0)=1$. This scalar product is positive:
$$(e_i,e_j)=0,\;i \ne j;\qquad (e_k,e_k)>0,\;k \in{\Bbb Z}.$$ The
irreducible $*$-representations of $U_q \frak{sl}_2$ in $V^{(l)}$, $l \in
{\cal L}_1,{\cal L}_2,{\cal L}_3$, are representations respectively from the
additional series, the principal unitary series, and the strange series (see
\cite{MMNNSU, VK}).

 It follows from theorem 3.9 that the selfadjoint linear operator
$\Box^{(0)}$ has a simple spectrum. Hence (see \cite{Y}), there exists a
Borel measure $dm(l)$ on the compact
$${\cal L}_0=\{l \in{\cal L}|\;c_1 \le-\frac{{\rm sh}(lh/2){\rm
sh}((l+1)h/2)}{{\rm sh}^2(h/2)}e^{h/2}\le c_2 \}$$
and a unitary operator $u:D(U)_q \to L^2(dm)$ such that $uf_0=1$, and for
all $f \in D(U)_q$
$$u \Box^{(0)}f=e^{h/2}\frac{{\rm sh}(lh/2){\rm sh}((l+1)h/2)}{{\rm
sh}^2(h/2)}uf.$$
(Note that $dm(l)$ and $u$ are uniquely determined by the properties listed
above.)

 Consider the completions $\overline{V^{(l)}}$ of the pre-Hilbert spaces
$V^{(l)}$, $l \in{\cal L}_0$, and the direct integral
$\overline{V}=\bigoplus \int \limits_{{\cal L}_0}V^{(l)}dm(l)$ (see
\cite{BR}). By the construction, the Casimir element $\Omega$ is well
defined in the Hilbert space $\overline{V}$, and $c_1 \le -q^{-1}\Omega \le
c_2$. Now it suffices to present an isometric operator $i:D(U)_q
\hookrightarrow \overline{V}$ with $i \Omega f=\Omega if$ for all $f \in
D(U)_q$.

 Let $l \in{\cal L}$, and $\Phi_l(x)$ a function on $q^{-2{\Bbb Z}}$ which
is a solution of the boundary problem
$$\Box^{(0)}\Phi_l(x)=e^{h/2}\frac{{\rm sh}(lh/2){\rm
sh}((l+1)h/2)}{{\rm sh}^2(h/2)}\Phi_l(x),\qquad \Phi_l(1)=1.\eqno(5.9)$$
In this case the distribution $\Phi_l(x)$, $x=(1-zz^*)^{-1}$, in the quantum
disc is a solution of the equation
$$\Omega \Phi_l=\frac{{\rm sh}(lh/2){\rm sh}((l+1)h/2)}{{\rm
sh}^2(h/2)}\Phi_l.$$

 The following lemma is a direct consequence of the definitions.

\medskip

\begin{lemma}\hfill
\begin{enumerate}
\item For each $l \in{\cal L}$ there exists a unique injective morphism of
$U_q \frak{sl}_2$-modules $i^{(l)}:V^{(l)}\to D(U)_q^\prime$ such that
$i^{(l)}e_0=\Phi_l$.
\item For each $l \in{\cal L}$ there exists a unique linear operator
$j^{(l)}:D(U)_q \to \overline{V}^{(l)}$ such that for all $f \in D(U)_q$, $v
\in \overline{V}^{(l)}$ one has $(j^{(l)}f,v)=\displaystyle \int
\limits_{U_q}(i^{(l)}v)^*fd \nu$.
\end{enumerate}
\end{lemma}

\medskip

\begin{proposition}\hfill
\begin{enumerate}
\item The linear operator $i:D(U)_q \to \bigoplus \int \limits_{{\cal
L}_0}\overline{V}^{(l)}dm(l)$, $i:f \mapsto j^{(l)}f$ is
isometric and $i \Omega f=\Omega if$ for all $f \in D(U)_q$.
\item The isometry $i$ is extendable by a continuity up to a unitary
operator $\overline{i}:L^2(d \nu)_q \to \bigoplus \int \limits_{{\cal
L}_0}\overline{V}^{(l)}dm(l)$.
\end{enumerate}
\end{proposition}

\smallskip

 {\bf Proof.} The relation $i \Omega f=\Omega if$, $f \in D(U)_q$, and the
invariance of the scalar product $[f_1,f_2]\stackrel{\rm
def}{=}(if_1,if_2)$, $f_1,f_2 \in D(U)_q$, follow directly from the fact that
$i^{(l)}$ are morphisms of $U_q \frak{sl}_2$-modules. Now one can deduce
from the uniqueness of the unitary operator $u$ that
$$uf(x)= q^2\int \limits_1^\infty
\overline{\Phi_l(x)}f(x)d_{q^{-2}}x.\eqno(5.10)$$
Hence $[f,f]=\int \limits_{U_q}f^*fd \nu$ for all $f \in D(U)_q$ which are
functions of $x=(1-zz^*)^{-1}$. The sesquilinear form $[f,f]-\int
\limits_{U_q}f^*fd \nu$ is zero on a subspace of finite functions of the
form $f(x)$. One can readily apply theorem 3.9 that such sesquilinear form
is identically zero. Thus the isometricity is proved. It remains to show
that the linear subspace $iD(U)_q$ is dense in the Hilbert space $\bigoplus
\int \limits_{{\cal L}_0}\overline{V}^{(l)}dm(l)$. Note that its closure $L$
is an invariant subspace of the multiplication operator by the function
$\frac{\textstyle{\rm sh}(lh/2){\rm sh}((l+1)h/2)}{\textstyle{\rm
sh}^2(h/2)}$. This function is real on the compact ${\cal L}_0$ and
separates its points. Hence $L$ is a common invariant subspace for
multiplication operators by functions $\varphi \in L^\infty(dm)$. It remains
to observe that the subspaces $j^{(l)}D(U)_q$ are dense in
$\overline{V}^{(l)}$ for all $l \in {\cal L}_0$ due to the injectivity of
the "conjugate" operators $i^{(l)}:V^{(l)}\hookrightarrow D(U)_q^\prime$. It
is well known (see \cite{BR}) that the listed properties of $L$ imply
$L=\bigoplus \int \limits_{{\cal L}_0}\overline{V}^{(l)}dm(l)$.  \hfill
$\Box$

\medskip

\begin{corollary} $c_1 \le-\Box \le c_2$.
\end{corollary}

\medskip

 This finishes the proof of \cite[lemma 3.1]{SSV1}.

\bigskip

\section*{Appendix. Covariant *-algebra ${\rm \bf Fun}\bf (U)_q$}

 Consider the vector space
$${\rm Fun}(U)_q \stackrel{def}{=}{\rm Pol}({\Bbb C})_q+D(U)_q \subset
D(U)_q^\prime.$$
It is easy to show that the multiplication, the action of the Hopf algebra
$U_q \frak{sl}_2$ and the involution $*$ can be transferred by a continuity
from ${\rm Pol}({\Bbb C})_q$ onto ${\rm Fun}(U)_q$. The elements $z,z^*,f_0$
generate this covariant *-algebra. The complete list of relations has the
form
$$z^*z=q^2zz^*+1-q^2;\quad z^*f_0=f_0z=0;\quad f_0\cdot f_0=f_0.$$
The action of the Hopf algebra $U_q \frak{sl}_2$ is determined by (1.1),
(2.4), (3.5), and the relation $H \cdot f_0=0$. The involution $*$ possesses
the following obvious properties: $*:z \mapsto z^*$, $*:f_0 \mapsto f_0$.

 One has the exact sequence of morphisms of covariant *-algebras
$$0 \to D(U)_q \to{\rm Fun}(U)_q \to{\rm Pol}({\Bbb C})_q \to 0.$$

 Every element $f \in D(U)_q^\prime$ admits a decomposition into series
$$f=\sum_{jk}a_{jk}(f)z^j \cdot f_0 \cdot z^{*k},\quad a_{jk}\in{\Bbb C}.$$
This decomposition is unique. $f \in D(U)_q^\prime$ is a finite function
iff ${\tt card}\{(j,k)\in{\Bbb Z}_+^2|\;a_{jk}(f)\ne 0 \}<\infty.$

\end{document}